\documentclass[sort&compress,preprint,12pt]{article}
\usepackage{authblk}
\usepackage{graphicx}
\usepackage[letterpaper,top=2cm,bottom=2cm,left=3cm,right=3cm,marginparwidth=1.75cm]{geometry}

\usepackage{amssymb}
\usepackage{amsmath}
\usepackage{color}
\usepackage{hyperref}

\usepackage{lineno}

\DeclareMathOperator{\diag}{diag}

\DeclareMathOperator{\sign}{sign}

\title{Quasiperiodic perturbations of Stokes waves: \\Secondary bifurcations and stability}

\author[1]{Sergey A. Dyachenko}

\affil[1]{ Department of Mathematics, SUNY at Buffalo, Buffalo, NY, USA}

\author[2]{Anastassiya Semenova}

\affil[2]{ Department of Applied Mathematics, University of Washington, Seattle, WA, USA}

\date{}
\begin{document}
\maketitle
\begin{abstract}
We develop a numerical method based on canonical conformal variables to study two eigenvalue problems for operators fundamental to finding a Stokes wave and its stability in a 2D ideal fluid with a free surface in infinite depth.
We determine the spectrum of the linearization operator of the quasiperiodic Babenko equation, and provide new results for eigenvalues and eigenvectors near the limiting Stokes wave identifying new bifurcation points via the Fourier-Floquet-Hill (FFH) method.
We conjecture that infinitely many secondary bifurcation points exist as the limiting Stokes wave is approached.
The eigenvalue problem for stability of Stokes waves is also considered.
The new technique is extended to allow finding of quasiperiodic eigenfunctions by introduction of FFH approach to the 
canonical conformal variables based method.
Our findings agree and extend existing results for the Benjamin-Feir, high-frequency and localized instabilities.
For both problems the numerical methods are based on Krylov subspaces and do not require forming of operator matrices. Application of each operator is pseudospectral employing the fast Fourier transform (FFT), thus enjoying the benefits of spectral accuracy and $O(N\log N)$ numerical complexity. 
Extension to nonuniform grid spacing is possible via introducing auxiliary conformal maps.
\end{abstract}




\section{Introduction}
Waves on the surface of the ocean appear due to multiple sources ranging from large scales such as seismic disturbances, motion of celestial bodies to the small scales dominated by forces of surface tension and dissipation.
Ocean and atmosphere are intertwined systems, and wind is one of the primary forces responsible for wave generation.
The details of generation of water waves by wind remain elusive, and it is generally thought to be a combination of the Miles-Phillips mechanism~\cite{miles1957generation,phillips1957generation} together with ideas of Jeffreys' sheltering theory~\cite{jeffreys1925formation}.
Nevertheless, ocean surface can be decoupled from atmospheric effects when long ocean waves are considered. 
As waves travel away from an epicenter of a storm, the long waves can be viewed as almost unidirectional periodic waves.  

Traveling waves propagating at a fixed speed without changing shape are called the Stokes waves. 
They are special solutions of $2$D Euler equations originally discovered by Stokes in~\cite{stokes1847theory,stokes1880} and later studied in many works~\cite{michell1893,nekrasov1921waves,schwartz1974computer,grant1973singularity,williams1981limiting, williams1985tables,longuet1978theory,longuet1977theory,cowley1999formation,longuet2008approximation}.
The Stokes waves are a one-parameter family of solutions, and steepness is commonly used for parameterization. It is denoted by $s = H/\Lambda$, and is defined as the ratio of crest-to-trough height $H$, over the wavelength $\Lambda$. The steepness ranges from zero to $s_{lim} = 0.14106348398\ldots$ (see also
Refs.~\cite{chandler1993computation,maklakov2002almost,dyachenko2014complex,lushnikov2017new,dyachenko2022almost}) in the limiting Stokes wave, or the wave of greatest height, which forms an angular crest with $120^\circ$.
Existence of such limiting wave has been proven in the works of~\cite{toland1978existence, amick1982stokes,plotnikov2002proof}.

The primary branch of Stokes waves contains solutions having one wavecrest per period. A secondary (double-period) bifurcation from the primary branch leads to a different branch of solutions having double the period of the original Stokes wave~\cite{chen1980numerical,zufiria1987non,wilkening2021spatially,wilkening2022spatially}. Each period contains a pair of wavecrests that are distinct.
Similarly, the secondary bifurcations happen for other ratios relative to primary period (e.g. tripling), and such bifurcations are effectively handled via the Floquet theory.

Dynamics of Stokes waves can be studied in the framework of Euler equations. 
The Hamiltonian framework to study the time-evolution of potential flow of an ideal fluid with free surface has been established in~\cite{zakharov1968stability}.
Equations governing the motion of a fluid surface in conformal variables were derived in~\cite{ovsyannikov1973dynamika,tanveer1991singularities,tanveer1993singularities,dyachenko1996dynamics,dyachenko1996nonlinear},
and wave breaking was numerically explored in~\cite{baker2011singularities,dyachenko2016whitecapping}. 
Instabilities of Stokes wave have been studied in many works~\cite{benjamin1967disintegration, bridgesmielke, lighthill1965contributions,whitham1967non,zakharov1968stability,ZakharovOstrovskyPhysD2009,berti2022full, creedon2022ahigh, nguyenstrauss}
including numerical studies~\cite{deconinck2011instability,dosaev2017simulation,murashige2020stability,deconinck2022instability}.
Stokes waves are unstable to long wave perturbations, a.k.a the {\em subharmonic} perturbations, such as the Benjamin-Feir (BF) or modulational instability, first reported in~\cite{benjamin1967disintegration, benjamin1967instability, whitham1967non} and high-frequency instability~\cite{deconinck2011instability,creedon2022high,creedon2022complete,creedon2022ahigh}.
Recent theoretical studies of Benjamin-Feir instabillity for small amplitude waves are described  in~\cite{creedon2022high, berti2022full, nguyenstrauss}.
The {\em superharmonic} instability  corresponds to growth of co-periodic perturbations
and has been described in~\cite{longuet1978instabilities} and found numerically in~\cite{tanaka1983stability, longuet1997crest} with almost limiting waves considered in~\cite{murashige2020stability,korotkevich2022superharmonic, deconinck2022instability}.

In the present work, we consider two eigenvalue problems connected to Stokes waves. 
In the first problem (see Section~\ref{section:EigenS1}), we developed a numerical method that is based on Krylov subspaces to find eigenvalues and eigenfunctions of the linearized Babenko equation~\cite{babenko1987some}. It is spectrally accurate, matrix-free and has $O(N\log N)$ numerical complexity, where $N$ is the number of Fourier modes approximating a Stokes wave. 
We extend the method with FFH approach~\cite{deconinck2006computing} to study quasiperiodic eigenfunctions.
Appearance of additional zero eigenvalue is found to occur at either a turning point of speed, or at a secondary bifurcation from the primary branch of Stokes waves.
We report new secondary bifurcations (see also~\cite{wilkening2022spatially}) and conjecture that infinite number of bifurcations may occur as the limiting Stokes wave is approached. Novel results are found for eigenfunctions of the linearized Babenko equation, and qualitative similarity with the classical Sturm-Louisville problem is observed. The eigenfunctions become strongly concentrated near the wavecrest, and may be related to the structure of the transition layer around the wavecrest~\cite{dyachenko2022almost}. The linearized Babenko operator is self-adjoint, and thus the set of its eigenfunctions is complete. The eigenfunctions of Babenko equation may be significant for theoretical study of wave-breaking phenomena similar to~\cite{dyachenko2013logarithmic,lushnikov2013beyond} to describe numerical observations in Refs.~\cite{dyachenko2016whitecapping,deconinck2022instability}.
In the second problem (see Section~\ref{section:eigenQEP}) we seek eigenvalues and eigenfunctions for the equations of motion linearized around the Stokes wave. 
We extended the canonical conformal variables based method (CCVM) in~\cite{dyachenko_semenova2022} with the Fourier-Floquet-Hill approach to study stability of Stokes waves to quasiperiodic disturbances. It is shown how a quadratic eigenvalue problem (QEP) with a pair of self-adjoint, and a skew-adjoint operators can be obtained. The QEP is solved efficiently via Krylov subspace based method avoiding matrix formation, having spectral accuracy and $O(N \log N)$ complexity. 
We provide examples that confirm and extend the stability results for
superharmonic instabilities that were found by linearization in Tanveer-Dyachenko variables~\cite{korotkevich2022superharmonic}, and provide new results for Benjamin-Feir, high-frequency and localized instabilities of Stokes waves. Additional results obtained by the present method are reported in the work~\cite{deconinck2022instability}.
%
%
 
%

\section{Euler equations with free surface}

\begin{figure}[htp]
\includegraphics[width=0.95\textwidth]{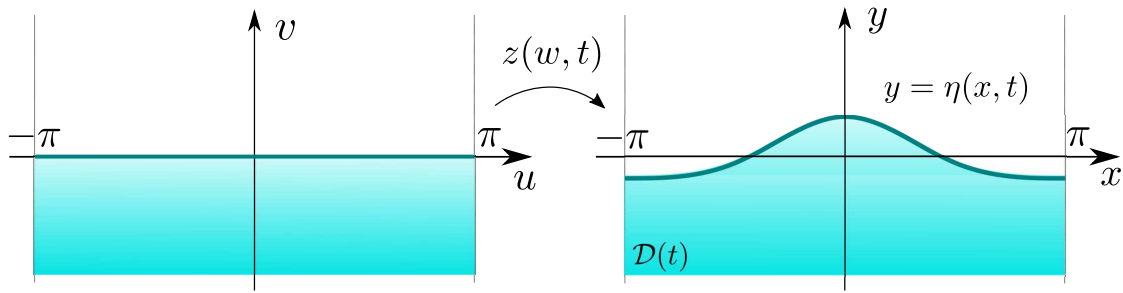}
\caption{The time-dependent conformal transformation $z(w,t) = x(w,t)+iy(w,t)$ maps $2\pi$-periodic semi-infinite strip in $w = u+iv$-plane (left panel) to $2\pi$-periodic infinitely deep fluid domain $\mathcal{D}(t)$ (right panel). 
The line $v = 0$ is mapped to the free surface $y = \eta(x,t)$.
}
\label{fig:confmap}
\end{figure}

We consider a potential flow of an ideal fluid of infinite depth with a free surface at $z = x + i\eta(x,t)$ (see Fig.~\ref{fig:confmap}). 
The fluid velocity is given by ${\bf v} = \nabla \varphi$ where $\varphi(x,y,t)$ is the velocity potential and $\nabla = \left(\partial_x, \partial_y\right)^T$.
By virtue of incompressibility $\varphi$ satisfies the Laplace equation $\Delta \varphi = 0$ in the fluid domain $\mathcal{D}(t)$
and the following boundary conditions are imposed at the free surface $y=\eta(x,t)$ and at the bottom of the fluid domain,
\begin{align}
&\dfrac{\partial \eta}{\partial t} + \left.\left(\dfrac{\partial \eta}{\partial x}\dfrac{\partial \varphi}{\partial x}  - \dfrac{\partial \varphi}{\partial y}\right)\right|_{y = \eta(x,t)} = 0 \label{bc1} \\
&\left.\left( \dfrac{\partial \varphi}{\partial t} + \dfrac{1}{2}\left(\nabla \varphi\right)^2 \right)\right|_{y = \eta(x,t)} +
g\eta = 0 \label{bc2} \\ 
&\left.\dfrac{\partial \varphi}{\partial y}\right|_{y\to -\infty} = 0 \label{bc3}
\end{align}
where $g$ is gravity acceleration. We introduce 
the surface potential, $\psi(x,t) = \varphi(x,\eta(x,t),t)$, and the variables $\psi(x,t)$ and $\eta(x,t)$
are the canonical variables as shown in~\cite{zakharov1968stability}.

\subsection{Hamiltonian and conformal variables}
The system~\eqref{bc1}--\eqref{bc3} together with Laplace equation is described by methods of Hamiltonian mechanics. 
The Hamiltonian is given by,
\begin{align}
    \mathcal{H} = \frac{1}{2} \int\limits_{-\pi}^\pi \int\limits_{-\infty}^{\eta(x,t)} \left(\nabla \varphi\right)^2\,dxdy + \frac{g}{2}\int\limits_{-\pi}^{\pi} \eta^2 \,dx,\label{HamilPhys}
\end{align}
which is the sum of kinetic and potential energy of the fluid.
Following the work~\cite{DyachenkoEtAl1996}, we express  Hamiltonian as surface integrals as follows:
\begin{align}
    \mathcal{H} = \frac{1}{2} \int\limits_{-\pi}^\pi 
    \psi \left.\frac{\partial \varphi}{\partial {\bf n}}\right|_{y = \eta(x,t)} \,dx + \frac{g}{2}\int\limits_{-\pi}^{\pi} \eta^2 \,dx,\label{HamPhys}
\end{align}
where ${\bf n}$ is the normal to the free surface, and Dirichlet-to-Neumann operator must be introduced to express the normal derivative of the potential at the free surface through $\psi$ and $\eta$. 
 The variables $\psi$ and $\eta$ are the canonical momentum and coordinate in physical variables, and the system:
\begin{align}
    \frac{\partial \eta}{\partial t} = \frac{\delta \mathcal{H}}{\delta \psi}, \quad 
    \frac{\partial \psi}{\partial t} = - \frac{\delta \mathcal{H}}{\delta \eta},
\end{align}
is equivalent to the system~\eqref{bc1}-\eqref{bc3}. Instead of writing the Hamiltonian in the the variables $\psi$,$\eta$ we follow~\cite{DyachenkoEtAl1996} and introduce a conformal mapping
from $w = u+iv\in\mathbb{C}^-$ to $\mathcal{D}(t)$ denoted by
$z(w,t) = x(w,t)+iy(w,t)$ (see Fig.~\ref{fig:confmap}).
In conformal variables, the free surface is described in parametric form,
\begin{align}
z(u,t) = x(u,t) + iy(u,t) = u - \hat H y(u,t) + iy(u,t),    
\end{align}
where $\hat H$ is the circular Hilbert transform defined as $\hat H e^{iku} = i\sign{(k)} e^{iku}$.

We abuse notation to define, $\psi(u,t) := \varphi(x(u,t),y(u,t),t)$ and $y(u,t):=\eta(x(u,t),t)$. The Hamiltonian in the conformal variables can be written explicitly as follows,
\begin{align}
    \mathcal{H} = \frac{1}{2}\int\limits_{-\pi}^{\pi} \psi \hat k \psi \,du + \frac{g}{2} \int\limits_{-\pi}^{\pi} y^2 x_u\,du,
\end{align}
with $\hat k = -\partial_u \hat H$, hence the Fourier symbol for $\hat k$ is $|k|$.
The Lagrangian in the conformal variables is given by,
\begin{align}
\mathcal{L} &= \int \psi \eta_t \,dx - \mathcal{H} = \int \psi \left(y_t x_u - y_u x_t\right)\,du - \mathcal{H}, \label{noncan}\\
\mathcal{L} &=
\int  \left(x_u\psi - \hat H\left[y_u\psi\right] \right)y_t \,du - \mathcal{H},
\end{align}
and canonical conformal momentum $\mathcal{P} := x_u\psi - \hat H[y_u \psi] =: \hat \Omega_{21}^\dagger \psi$ is introduced as in~\cite{DLZ1995Fivewave}, and write:
\begin{align}
    \mathcal{L} = \int \mathcal{P} y_t\,du - \mathcal{H}(\mathcal{P},y), \,\,\,\mathcal{H} = \frac{1}{2}\int \mathcal{P} \hat R_{12}\hat k \hat R^\dagger_{12} \mathcal{P}\,du + \frac{g}{2}\int y^2x_u \,du, 
\end{align}
with $\psi = \frac{x_u \mathcal{P} + \hat H \left[y_u\mathcal{P} \right]}{|z_u|^2} = :\hat R_{12}^\dagger \mathcal{P}$. Here the definitions of operators $\hat \Omega_{21}$ and $\hat R_{12}$ are in agreement with~\cite{LushnikovEtAl2019} and $\hat R_{12}\hat\Omega_{21} = 1$, with dagger denoting the adjoint operator.
In the variables $\mathcal{P}$ and $y$, the
system acquires canonical structure,
\begin{align}
    \frac{\partial y}{\partial t} = \frac{\delta H}{\delta \mathcal{P}}, \quad
    \frac{\partial \mathcal{P}}{\partial t} = - \frac{\delta H}{\delta y}.
\end{align}
The canonical variables $\mathcal{P}$ and $y$ naturally appear in  
numerical method for stability problem which is described in the Section~\ref{section:eigenQEP}.

\section{Implicit form and Babenko equation}

Equations of motion in implicit form are determined from the stationary action principle with Lagranian defined in~\eqref{noncan},
\begin{align}
\mathcal{S} &= \int \left(\mathcal{L} + \mathcal{C}\right) dt, \\
\mathcal{C} &= \int f(y - \hat H [x-u])\,du,
\end{align} 
where $\mathcal{S}$ is action, $\mathcal{C}$ is a Lagrangian constraint, and $f$ is the Lagrange multiplier ensuring that $x$ and $y$ are related by the Hilbert transform.
The equations of motion are 
obtained from extremizing action which implies that,
\begin{align}
\frac{\delta \mathcal{S}}{\delta x} = 0, \,\, \frac{\delta \mathcal{S}}{\delta y} = 0, \,\, \frac{\delta \mathcal{S}}{\delta \psi} = 0 \,\,\,\,\mbox{and}\,\,\,\,\frac{\delta \mathcal{S}}{\delta f} = 0
\end{align}
and solving for the Lagrange multiplier $f$ (see also Ref.~\cite{DyachenkoEtAl1996}) yields,
\begin{align}
&y_t x_u - y_u x_t = -\hat H\psi_u, \label{baseKC}  \\
&\psi_t x_u - x_t \psi_u + \hat H\left[-y_u\psi_t + y_t \psi_u \right] +
g\left[yx_u - \hat H\left(yy_u\right)\right] = 0. \label{baseDC}
\end{align} 
We switch to the reference frame traveling with the speed of a Stokes wave by considering the following ansatz,
\begin{align}
y(u,t) &\to y(u-ct,t), \quad \mbox{and} \quad x(u,t) \to u - \hat H y(u-ct,t), \label{ans1}\\ 
\psi(u,t) &\to\phi(u-ct, t) -  c\hat H y(u-ct,t) , \label{ans2}
\end{align}
where $c$ is the velocity of the moving frame. 

Substitution of~\eqref{ans1}-\eqref{ans2} into the equations of motion~\eqref{baseKC}--\eqref{baseDC} leads to the following
system in the moving reference frame:
\begin{align}
&x_u y_t - y_u x_t =  - \hat H\left(\phi_u - c \right), \label{KCm}\\ 
&x_u \phi_t - x_t\phi_u - \hat H\left[y_u \phi_t - y_t \phi_u\right]
+ 2cx_t= c^2\hat k y - g\left(x_uy - \hat H \left[yy_u\right]\right) .\label{DCm}
\end{align}
The Stokes wave is a stationary solution of the latter system. It is found
by setting $x_t = y_t = \phi_t = 0$ in\eqref{KCm}-\eqref{DCm} which leads to the relations:
\begin{align}
&\phi_u = 0 ,\quad\mbox{or}\quad \psi_u = c(x_u - 1) = -c\hat Hy_u = c\hat k y, \\
&\hat S y := c^2\hat k y - g\left[x_uy - \hat H (yy_u)\right] = 0. \label{Stokes}
\end{align}
The latter equation is also sometimes referred to as the Babenko equation~\cite{babenko1987some}.  
The operator $\hat S(y)$ is linearized around a Stokes wave $y(u)$ to form $\hat S_1(y)\delta y$ as follows,
\begin{align}
\hat S_1(y)\delta y := \left(c^2\hat k - g\right)\delta y - 
g\left[y\hat k \delta y + \delta y \hat k y + \hat k (y\delta y) \right]\label{S1}
\end{align}
which is {\em hermitian}, i.e. $\hat S_1^\dagger = \hat S_1$ with standard inner product. 

The discussion of eigenvalue problem for $\hat S_1$, its significance as well as new results and the numerical method is presented in the Section~\ref{section:EigenS1} of this paper.

\section{Linearization of the implicit equations of motion}

Linear stability is determined by the eigenvalue spectrum of linearization operator of the system~\eqref{KCm}--\eqref{DCm} around a Stokes wave, $y(u)$.
After substitution of the ansatz, 
$y(u,t) \to y(u) + \delta y(u,t)$ and $\phi(u,t) \to 0 + \delta \phi(u,t)$
in~\eqref{KCm}-\eqref{DCm} and keeping only the linear terms in $\delta y$ and $\delta \phi$ we obtain the linearized equations of motion, 
\begin{align}
&x_u \delta \phi_t - \hat H\left[ y_u\delta\phi_t \right] - 2c\hat H \delta y_t = \hat S_1(y)\delta y, \label{linDCm} \\
&x_u \delta y_t + y_u \hat H\delta y_t = -\hat H\delta\phi_u.\label{linKCm}
\end{align}
The equations~\eqref{linDCm}--\eqref{linKCm} are convenient to express in the operator matrix form:
\begin{align}
&\hat L
\left[
\begin{array}{c}
\delta \phi_t \\
\delta y_t
\end{array}
\right] =
\hat M
\left[
\begin{array}{c}
\delta \phi \\
\delta y
\end{array}
\right], \label{linsys1}
\end{align}
with the following operator matrices:
\begin{align}
\hat M =
\left[
\begin{array}{c|c}
\hat k  & 0 \\ \hline
 0  & \hat S_1
\end{array}
\right],\,\, \hat L =
\left[
\begin{array}{c|c}
0  & \hat \Omega_{21} \\ \hline
\hat \Omega_{21}^\dagger  & -2c\hat H
\end{array}
\right].
\label{OperLin}
\end{align}
Equivalently, the first-order system~\eqref{linsys1} given by
\begin{align}
&\hat \Omega_{21}^\dagger \delta \phi_t - 2c\hat H \delta y_t - \hat S_1 \delta y = 0, \label{dcl}\\
&\hat \Omega_{21}         \delta y_t = \hat k \delta \phi\label{kcl}
\end{align}
is expressed as a second order equation by solving~\eqref{kcl} for $\delta \phi$ and substituting the result into~\eqref{dcl},
\begin{align}
\hat \Omega_{21}^\dagger\hat k^{-1} \hat\Omega_{21} \delta y_{tt} - 2c\hat H \delta y_t - \hat S_1 \delta y = 0. \label{QEP}
\end{align}
The equation~\eqref{QEP} leads to a quadratic eigenvalue problem (QEP) with substitution $\partial_t \to \lambda$ which determines the stability of Stokes waves.

\subsection{Dispersion of linear gravity waves}

The linear waves are stable, and the dispersion relation of linear waves in laboratory frame 
is given by:
\begin{align}
\omega(k) = \pm\sqrt{g|k|}.
\end{align}
This relation is trivially found from the QEP by noting 
that for a flat surface:
\begin{align}
    &\hat \Omega_{21} = \hat \Omega_{21}^\dagger = 1, \\
    &\hat S_1 = c^2 \hat k - g.
\end{align}
Then, the substitution of the plane wave in the form $\delta y \sim e^{i(ku - \omega t)}$ leads to 
\begin{align}
    \omega(k) = ck \pm \sqrt{g|k|},
\end{align}
where the term $ck$ is a Doppler frequency shift from moving to the traveling frame of reference.

We discuss the details of the numerical method for the system~\eqref{dcl}--\eqref{kcl} and new results obtained for the stability problem of large amplitude waves in the section~\ref{section:eigenQEP}.

\section{Fourier-Floquet-Hill method in conformal variables}

In order to consider quasiperiodic 
eigenfunctions we use the Fourier-Floquet-Hill method (FFH) described in~\cite{deconinck2006computing}
and implemented in ~\cite{deconinck2011instability,murashige2020stability}. We seek quasiperiodic eigenfunctions  $f(u)$ described by the trigonometric series:
\begin{align}
    f(u) = \tilde{f} (u) e^{i\mu u} = \sum\limits_{k=-\infty}^{\infty} \tilde{f}_k e^{i\left(k+\mu\right)u},
\end{align}
where $\tilde{f}(u)$ is a $2\pi$-periodic function.
We consider the operators of linearized equations of motion when applied to a quasiperiodic function $f$,
\begin{align}
    \hat k f = \hat k \left(\sum\limits_{k=-\infty}^{\infty} \tilde{f}_k e^{i\left(k+\mu\right)u} \right) = 
    e^{i\mu u}\sum\limits_{k=-\infty}^{\infty}  |k + \mu| \tilde{f}_k  e^{iku} = e^{i\mu u} \hat k_\mu \tilde{f}.
\end{align}
Hence, the operator $\hat k$
acts as multiplier $|k+\mu|$ in the Fourier space.
The linearized Babenko operator applied to $f$ is given by,
\begin{align}
    \hat S_1 f = e^{i\mu u} S_{1,\mu} \tilde{f} = 
    e^{i\mu u }\left[\left(c^2\hat k_\mu - g\right) \tilde{f} - g\left( 
    y\hat k_{\mu} \tilde{f} + \tilde{f}\hat k y + \hat k_\mu\left(y\tilde{f}\right)\right) \right],
\end{align}
and finally:
\begin{align}
    &\hat \Omega_{21} f = x_u f + y_u\hat H f = e^{i\mu u}
    \left[x_u \tilde{f} + y_u \hat H_\mu \tilde{f} \right] = e^{i\mu u} \hat \Omega_{21,\mu} \tilde{f}, \\
    &\hat \Omega^\dagger_{21} f = x_u f - \hat H \left(y_u f\right) = e^{i\mu u}
    \left[ x_u \tilde{f} - \hat H_\mu\left(y_u \tilde{f} \right)\right] = e^{i\mu u} \hat \Omega_{21,\mu}^\dagger \tilde{f}, 
\end{align}
where $\hat H_{k,\mu} = i\sign{(k+\mu)}$ is the Fourier multiplier for the quasiperiodic Hilbert transform, $\hat H_\mu$.

\section{\label{section:EigenS1}Results: Eigenvalues of $\hat S_1$.}

\begin{figure}[htp]
\includegraphics[width=0.495\textwidth]{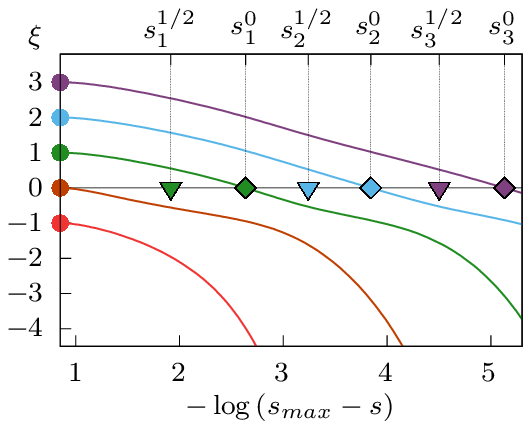}
\includegraphics[width=0.495\textwidth]{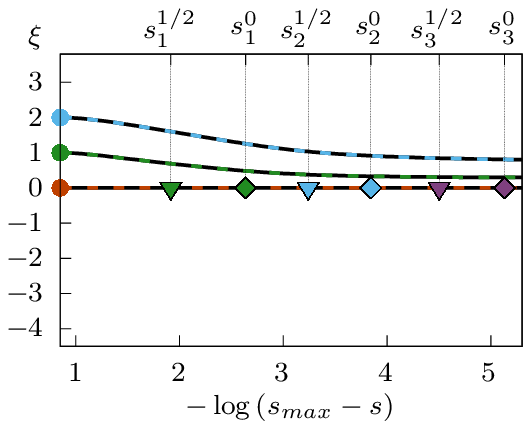}
\caption{The first few eigenvalues $\xi$ of $\hat S_1$ as a function of the Stokes wave steepness. At flat surface each eigenvalue has multiplicity $2$ (orange, green, blue and violet circles), except $-1$ which is simple (red circle). The double eigenvalues break into simple eigenvalues in finite amplitude Stokes waves.
(Left Panel) The eigenvalues originating from a flat surface at $\xi = -1,0,1,2,3$ are colored in red, orange, green, blue and violet respectively; and (Right Panel) 
the second eigenvalue originating from $\xi = 0,1,2$ are marked dashed red, dashed green and dashed blue respectively. The green, blue and violet diamonds correspond to bifurcation 
points at the turning points of speed, and the triangles denote the double-period bifurcations.
}
\label{fig:eigS1}
\end{figure}

Occurrence of zero eigenvalues in $\hat S_1$ (aside from the trivial zero associated with Galilean invariance) indicates appearance of extra eigenfunction in the null space of $\hat S_1$ and 
corresponds to a turning point of speed, or a secondary bifurcation from a branch of Stokes waves. We study the eigenvalues 
Since $\hat S_1$ is self-adjoint its eigenvalues are real, and linearization at a flat surface indicates that Fourier modes are the eigenfunctions as evident from,
\begin{align}
    \hat S_1 f \to \left(c^2\hat k - g\right) f.
\end{align}
Hence the eigenvalues are the integers $\xi = -1, 0 , 1, \ldots$ (with $c = 1$ and $g = 1$). Each eigenvalue except $-1$ has multiplicity $2$. The double eigenvalue $0$ (for $s=0$) signifies a primary bifurcation from the flat surface to a nonlinear Stokes wave (see Fig.~\ref{fig:eigS1}). In case of a flat surface
$\xi_k = |k-1|$ is an eigenvalue for $k = 0, 1, \ldots$. Each eigenvalue except $\xi_0 = -1$ has a $2$-dimensional eigenspace spanned by the functions $e^{\pm iku}$, or equivalently $\cos ku$ and $\sin ku$, associated with $\xi_k$.

Away from the flat surface, the eigenvalues of operator $\hat S_1$ are found numerically. Each eigenvalue becomes simple, except for the special cases when collisions of eigenvalues occur. A pair of eigenvalues originating from a double eigenvalue in a flat surface follow distinct paths: one eigenvalue corresponding to an even eigenfunction crosses the zero axis, whereas the other one remains positive. The Fig.~\ref{fig:eigS1} illustrates  typical behaviour of eigenvalues associated with even and odd eigenfunctions as a function of steepness.

For {\em regular waves}, the linearization operator $\hat S_1$ has a simple zero eigenvalue and associated eigenfunction $y_u$:
\begin{align}
    \hat S_1 y_u = \partial_u \left( \hat S y \right) = 0,
\end{align}
and corresponds to translational invariance of solutions to Babenko equation. 
%

At the bifurcation points (see Fig.~\ref{fig:eigS1}) the zero eigenvalues has double multiplicity, and the additional eigenpair is associated  
with either the turning points in the $c(s)$ graph,
or the secondary bifurcations from the primary branch (when subharmonics are considered). See also the table~\ref{table1} in the section~\ref{section:Bifurcation} for the few first secondary bifurcations and turning points.

The Fig.~\ref{fig:eigS1} illustrates the dependence of the first few eigenvalues $\xi_j$ of $\hat S_1$ on steepness, $s$. All the eigenvalues we found do not increase with steepness, i.e. $\partial_s \xi_j \leq 0$ for all $j$. Eventually half of eigenvalues cross the horizontal axis at turning points of speed at $s = s^{0}_{1}, s^{0}_{2}, \ldots$ and become negative. Since the number 
of eigenvalues is infinite, it is conjectured that 
zero eigenvalue appears infinitely many times as the 
extreme wave is approached, and is in agreement with~\cite{longuet1978theory}. 
Tracking the eigenvalues crossing the origin, one may find secondary bifurcations corresponding to 
appearance of irregular waves discovered in~\cite{chen1980numerical} and for finite depth in~\cite{vanden1983some}. In such case, the eigenvalue problem comes from linearization around  two periods of a Stokes wave.

\begin{figure}
    \centering
    \includegraphics[width=0.48\textwidth]{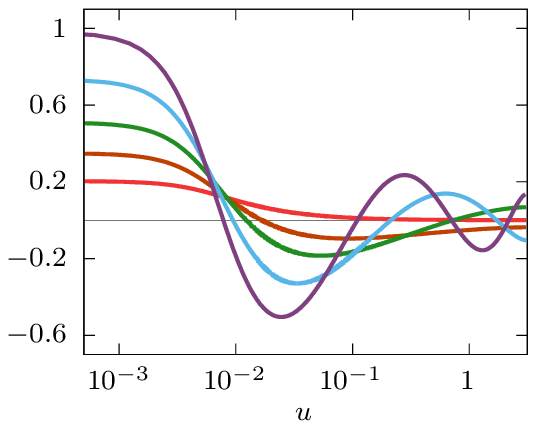}
    \includegraphics[width=0.48\textwidth]{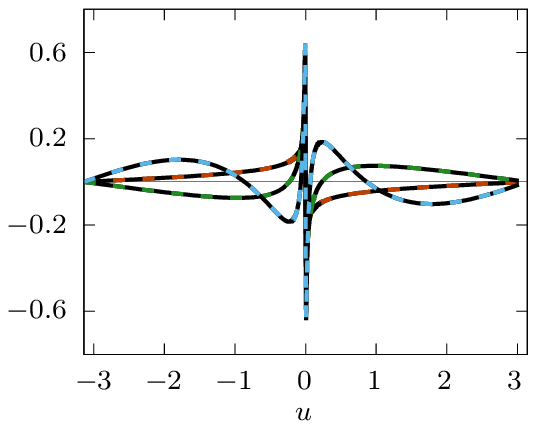}
    \caption{(Left Panel) The first five even eigenfunctions associated with eigenvalues at $s = 0.1388325$ (same colormap as in the left panel Fig.~\ref{fig:eigS1}), and (Right Panel) the first three odd eigenfunctions corresponding to the eigenvalues of the right panel of Fig.~\ref{fig:eigS1}.}
    \label{fig:eigfS1}
\end{figure}

Fig.~\ref{fig:eigfS1} shows the first few eigenfunctions for a large amplitude wave with steepness $s = 0.1388325$. In the case of a flat surface the even eigenfunctions are given by 
\begin{align}
    f_1(u) = 1, \quad f_2(u) = \cos{u}, \quad f_3(u) = \cos{2u}, \,\,\ldots.
\end{align}
The eigenfunctions continuously deform as we follow the primary branch and approach the limiting wave. Odd eigenfunctions continuosly deform from $\sin(ku)$, but the corresponding eigenvalues never cross zero. 

In addition we find that as the limiting wave is approached, the eigenfunctions become 
concentrated near the wave crest (see the Fig.~\ref{fig:eigfS1}) suggesting self-similar behaviour of the boundary 
layer near wave-crest discussed in~\cite{chandler1993computation,dyachenko2022almost}. In the remainder of the interval the eigenfunctions look qualitatively similar to the Fourier cosine/sine of the appropriate wavenumber.

\subsection{\label{section:Bifurcation}Bifurcation points of the linearized Babenko equation}
In this section we demonstrate several waves at bifurcations associated with turning points of speed ($\mu = 0$ coperiodic) and double period bifurcations ($\mu = 1/2$) in the primary branch. The result for double-period bifurcation point at $s^{\mu = 1/2}_1$ has been determined in~\cite{chen1980numerical,longuet1985bifurcation}, see also~\cite{zufiria1987non} for more bifurcation points corresponding to tripling and higher wavelength multiples.
We report a new bifurcation point at $s^{\mu = 1/2}_{3} = 0.141032049$ located between the first minimum and the second maximum of the Stokes wave speed.

\begin{table}[h!]
\begin{tabular}{l|l|l}
     \textbf{Notation} & \textbf{Steepness}, $s$ & \textbf{Description} \\ \hline  
       & $0$ & primary bifurcation \\
      $s^{\mu=1/2}_1$ & $0.128903\pm 2\times10^{-7}$ & double period bifurcation point \\
      $s^{\mu =0}_1$  & $0.138753\pm 4\times10^{-7}$ & coperiodic, turning point of speed \\
      $s^{\mu = 1/2}_2$ & $0.140487 \pm 9\times10^{-7}$ & double period bifurcation point \\
      $s^{\mu = 0}_2$ & $0.140920 \pm 1\times10^{-6}$ & coperiodic, turning point of speed \\
      $s^{\mu = 1/2}_3$ & $0.141032049\pm 5\times 10^{-9}$ & double period bifurcation point \\
      $s^{\mu = 0}_3$ &  $0.141056063 \pm 5\times 10^{-9}$ & coperiodic, turning point of speed
\end{tabular}
\caption{Steepness of waves at the first few turning points of speed and at the double-period bifurcations.}
\label{table1}
\end{table}
\noindent The Table~\ref{table1} indicates that there is a double-period bifurcation point enclosed between each turning point 
of speed,
\begin{align*}
0 < s^{\mu=1/2}_1 < s^{\mu = 0}_1 \,\,\,\mbox{and}\,\,\, s^{\mu=0}_1 < s^{\mu = 1/2}_2 < s^{\mu = 0}_2 \,\,\,\mbox{and}\,\,\,
s^{\mu=0}_2 < s^{\mu = 1/2}_3 < s^{\mu=0}_3.
\end{align*}
We conjecture that the same pattern extends to the limiting wave, and 
\begin{align}
    s^{\mu = 0}_{k} < s^{\mu = 1/2}_{k+1} < s^{\mu = 0}_{k+1}\quad\mbox{for all $k= 0,1,2,\ldots$}
\end{align}
where we defined $s^{\mu=0}_0 := 0$. Analogous result can be obtained for other values of $\mu$.

We would also like to refer the reader to the recent work~\cite{wilkening2022spatially} that uses singular value 
decomposition (SVD) of the Jacobian to identify bifurcation points. In addition, the latter work offers a way to determine 
fully nonlinear quasiperiodic traveling waves with two quasiperiods, and allow finding of secondary bifurcations from 
these quasiperiodic branches. Significance of subharmonic bifurcation points was emphasized in~\cite{vanden2014periodic} as 
a possible route to determining of solitary waves.
 
\subsection{Numerical Method}

We seek solution of the eigenvalue problem:
\begin{align}
\hat S_1 f = \xi f,
\end{align}
where $\xi$ is a real eigenvalue, and $f(u)$ is a 
$2\pi$-periodic eigenfunction.
The linear operator $\hat S_1$ is self-adjoint,
\begin{align}
\left(f\cdot\hat S_1 g\right) = \int f \hat S_1 g \,du = \int g\hat S_1 f  \,du =
\left(\hat S_1f\cdot g\right).
\end{align}

%
%
The eigenvalue problem can be written in terms of the auxiliary variable, $q$, using the auxiliary conformal map~\cite{lushnikov2017new}:
\begin{align}
    \tan\frac{u}{2} = L\tan\frac{q}{2}\,\,\,\mbox{and}\,\,\, 
    u_q = \frac{2L}{1+L^2 - \left(1-L^2\right)\cos{q}} = \frac{1}{q_u},
\end{align}
where $0<L\leq 1$ is a parameter controlling the distribution of grid points, and $u_q$ and $q_u$ are $2\pi$-periodic functions defined at  collocation points.
In the $q$ variable, the operator $\hat S_1$ is given by:
\begin{align}
    \hat S_1 f(q) = q_u\left[ \left(c^2\hat k_q - gu_q\right)f
    - g \left( y \hat k_q f + f \hat k_q y + \hat k_q \left(yf\right)\right) \right],
\end{align}
where $\hat k_q$ has the Fourier multiplier $|k|$ in the $q$-variable, i.e. $\hat k_q\left( e^{ikq}\right) = |k|e^{ikq}$.
It is convenient to introduce a self-adjoint operator $\hat A$,
\begin{align}
    \hat A f(q) := \left(c^2\hat k_q - gu_q\right)f
    - g \left( y \hat k_q f + f \hat k_q y + \hat k_q \left(yf\right)\right),
\end{align}
so the eigenvalue problem may be written as follows:
\begin{align}
\left[\hat A - \xi u_q \right]f(q) = 0, \quad\mbox{or}\quad
\left[q_u^{1/2} \hat A q_u^{1/2} - \xi \right]  h(q) = 0 \label{dirEig1}
\end{align}
here $h(q) = u_q^{1/2} f(q)$. The first form is a generalized eigenvalue problem with 
a hermitian pair of operators ($u_q$ is symmetric positive definite), or a regular eigenvalue problem with hermitian operator $q_u^{1/2} \hat A q_u^{1/2}$. The shift--inverse method (see also Ref.~\cite{saad1992numerical}) allows to determine the eigenvalue nearest to an arbitrary point $\sigma$
in the spectral plane. It amounts to solving a sequence of linear problems given by
\begin{align}
\left[q_u^{1/2} \hat A q_u^{1/2} - (\xi + \sigma) \right] f^{(n+1)} = f^{(n)}, \label{shiftinvertS1}
\end{align}
where $\sigma$ is the shift, and superscript denotes the iteration number.
The linear system~\eqref{shiftinvertS1} is solved 
by means of the minimal residual method (MINRES), and defines,
\begin{align}
    f^{(n+1)} = B_{\sigma} f^{(n)}\quad\mbox{where}\quad B_{\sigma} =\left[q_u^{1/2} \hat A q_u^{1/2} - (\xi + \sigma) \right]^{-1}.
\end{align}
The result, $f^{(n+1)}$ is a new basis vector for a subsequent Krylov subspace $\mathcal{K}_n (B_{\sigma}, f^{(0)}) = \mathrm{span}\{f^{(0)}, B_{\sigma} f^{(0)}, \dots, B^{n-1}_{\sigma} f^{(0)}\}.$
In the limit $n\to\infty$, the 
sequence of $f^{(n)}$ converges to the desired eigenfunction. The choice of the initial guess $f^{(0)}$ and the shift $\sigma$ is made via 
the continuation method from the linear waves.

Advanced numerical methods (see Refs.~\cite{lehoucq1989arpack,lehoucq1996arpack,stewart2002krylov}) are also available to seek multiple eigenvalues simultaneously. Note, that the operator matrix for $B_{\sigma}$ is never formed, and applying the linear operator in~\eqref{shiftinvertS1} to an arbitrary $2\pi$-periodic function via FFT is spectrally accurate and requires $O(m N\log N)$ flops, where $m$ is 
the number of MINRES iterations that it takes to solve the associated linear system to a given tolerance.
We considered up to $4\times1024^2$ Fourier modes
and uniform grid to seek double period eigenfunctions using FFH method, and up to $16384$ Fourier modes with nonuniform grid with $L \sim 0.005$.

\section{\label{section:eigenQEP}Results: Stability problem}

In this section we illustrate the FFH method with canonical conformal variables based approach~\cite{dyachenko_semenova2022} in a sequence of simulations and determine the stability spectrum of Stokes waves for quasi-periodic perturbations. The Stokes waves themselves are computed by means of the Newton method coupled with the conjugate residual method~\cite{yang2009newton,yang2010nonlinear}. The library containing the Stokes waves is available at~\texttt{stokeswave.org} and 
a discussion on how to compute Stokes waves can be found in  Refs.~\cite{dyachenko2014complex,lushnikov2017new,dyachenko2022almost}.

\begin{figure}
    \centering
    \includegraphics[width=0.43\textwidth]{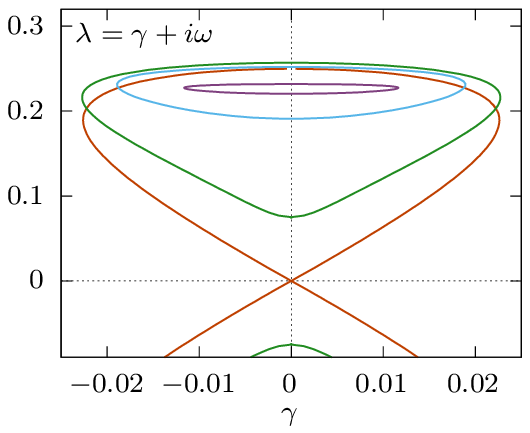}
    \includegraphics[width=0.56\textwidth]{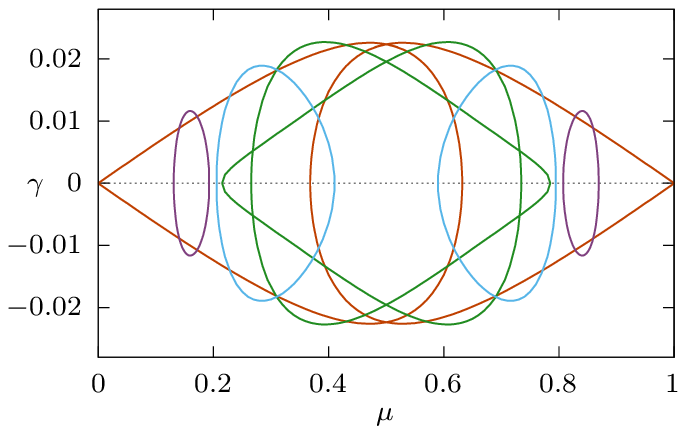}
    \caption{(Left Panel) The spectral plane showing BF instability in four Stokes waves with steepnesses
    $s = 0.095493$ (red), $s = 0.111408$ (green), $s = 0.119366$ (blue), and $s = 0.127324$ (violet). (Right Panel) The values of Floquet parameter $\mu$ vs real parts of eigenvalues $\gamma$ for the BF instabilities only. 
    }
    \label{fig:fig_eight}
\end{figure}

\subsection{Benjamin-Feir Instability}
The Benjamin-Feir (BF), or modulational instability is dominant for Stokes waves with steepness less than $s \approx 0.12894$. It has been well-studied, and has 
a rigourous theoretical description in small amplitude waves regime~\cite{nguyenstrauss}.
In~\cite{berti2022full,creedon2022ahigh} an asymptotic theory is developed giving excellent quantitative predictions for up to $s \approx 0.32$ ($ka \approx 0.1$). 
However, the BF instability for larger amplitude waves becomes challenging to study numerically since it is only present in a narrow band of the Floquet parameter $\mu$.
For example, in the work~\cite{dyachenko_semenova2022} $n=64$ subharmonics were considered to study Benjamin-Feir instability in Stokes waves up to $s = 0.12732395$ ($ka = 0.4$). It was found that no Benjamin-Feir instability is present for $s=0.12732395$ with $n=64$ subharmonics.
However, the FFH approach allows to scan arbitrary fine bands in Floquet parameter, and we demonstrate that BF instability persists for the wave with $s = 0.12732395$ albeit in a very narrow band in $\mu$ (see the violet curve in the Fig.~\ref{fig:fig_eight}).

In the left panel of Fig.~\ref{fig:fig_eight} we show the spectrum of BF instability as steepness of Stokes wave increases, and the range of $\mu$ values is presented in the right panel. 
We find that the remnant of BF instability becomes qualitatively similar to the high-frequency instability as the figure--eight splits and shrinks.

\begin{figure}
    \centering
    \includegraphics[width=0.506\textwidth]{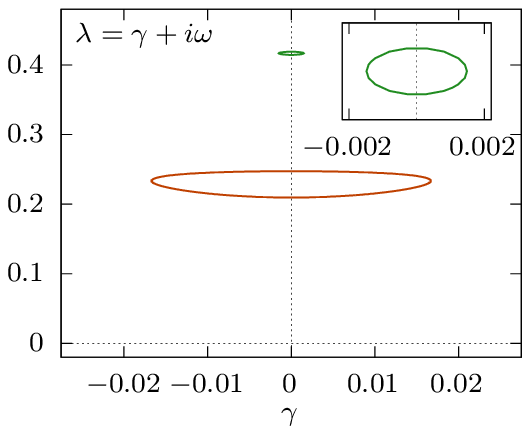}
    \includegraphics[width=0.477\textwidth]{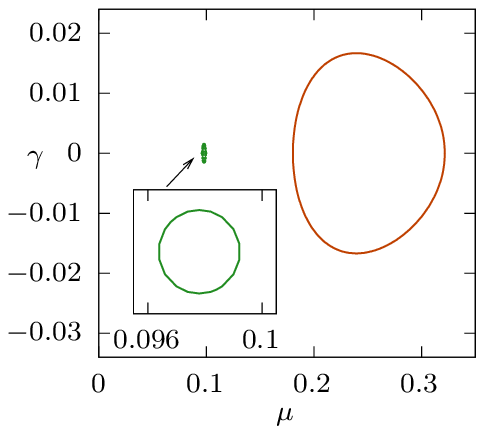}
    \caption{(Left Panel) The eigenvalues of~\eqref{quasiQEP} for the remnant of BF instability (orange) and the high-frequency instability (green) in a Stokes wave with $s = 0.1222625$, and (Right Panel) the growth rate, $\gamma$ as a function of the Floquet parameter, $\mu$, high-frequency instability is contained in the range $\mu\in(0.964,0.992)$.
    }
    \label{fig:hf1}
\end{figure}

\subsection{High-Frequency Instability}
The high-frequency instability discovered in~\cite{deconinck2011instability} is nontrivial 
to find because of the very limited range of the Floquet parameter, $\mu$, for which it is present.
In the left panel of Fig.~\ref{fig:hf1} we show the spectral plane ($\lambda = \gamma + i\omega$) for $s=0.12222625$ with the BF (red points) and  the high-frequency instability (green points). The range of $\mu$ values vs growth rate, $\gamma$, of eigenvalues for BF and high-frequency instabilities is shown in the right panel. 
When compared to the remnant of the BF instability in Fig.~\ref{fig:hf1}, the range of Floquet multipliers for resolving the high-frequency instability is even finer. For example, for a Stokes wave with $s = 0.1222625$ the entire instability bubble is contained within the interval $\mu \in (0.964,0.992)$ and would take a prohibitively large number of subharmonics unless FFH method is used. We observe that the growth rate for high-frequency instability can be significant, e.g 
for $s=0.122262$ it reaches up to $\gamma = 1.48\times10^{-3}$ (see the Fig.~\ref{fig:hf1}).

\subsection{Localized Instability}
\begin{figure}
    \centering
    \includegraphics[width=0.995\textwidth]{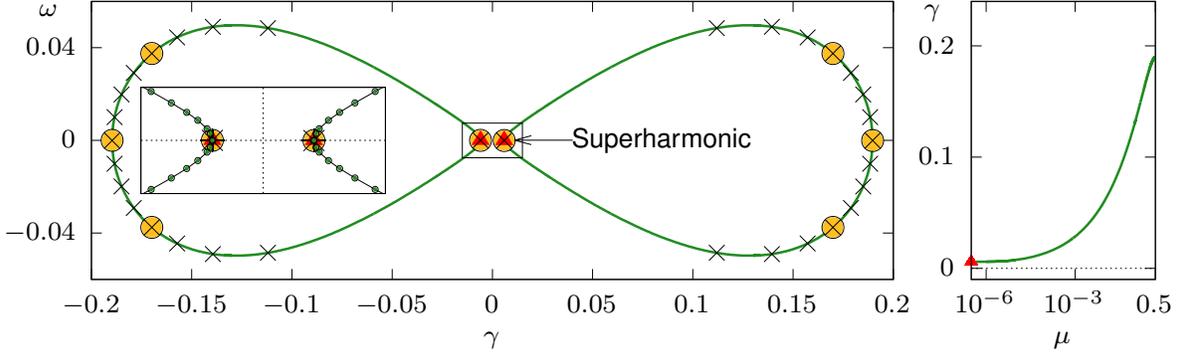}
    \caption{
    (Left Panel) Emergence of the localized instability in~\cite{deconinck2022instability} for Stokes wave with steepness $s = 0.13660515$. Eigenvalues presented are found by three methods: the canonical conformal variables-based method (CCVM) with FFH (solid green, green circles in inset); CCVM from~\cite{dyachenko_semenova2022} with $1$ (red triangles), $4$ (gold circles) and $16$ subharmonics (black crosses); matrix-forming linearization from~\cite{korotkevich2022superharmonic} (black plus signs).
    (Right Panel) Growth rate as a function of $\mu$ for $\mu \in [0,0.5]$ only positive values are shown. 
    }
    \label{fig:infty}
\end{figure}

The localized instability branch appears before the first extremum of the Hamiltonian, and fully detaches at its extremum~\cite{deconinck2022instability}. The localized branch dominates 
dynamics of steep waves, and induces wave-breaking in the nonlinear stage of time evolution of Stokes waves~\cite{dyachenko2016whitecapping,dyachenko_semenova2022}.
In Fig.~\ref{fig:infty} we explore the spectrum of dynamical system for the Stokes wave with $s = 0.1366051511$.
We found that the FFH method in canonical conformal variables is capable of extending the results of the numerical simulations in~\cite{dyachenko_semenova2022,korotkevich2022superharmonic}. It allows us to recover a continuous curve of eigenvalues, $\lambda(\mu)$, that corresponds to the localized instability branch~\cite{deconinck2022instability}.

We observe that the Floquet parameter has to be varied on a logarithmic scale around $\mu = 0$, so that the points on the $\lambda(\mu)$ curve are uniformly resolved in the vicinity of the point of superharmonic instability. We find that the results obtained with FFH approach together with canonical conformal variables based method are in agreement with the findings in~\cite{dyachenko_semenova2022,korotkevich2022superharmonic}.

\subsection{Numerical Method}
The QEP is obtained by substitution $\partial_t \to \lambda$ into~\eqref{QEP}, and is given by:
\begin{align}
    \left[ \lambda^2 \hat \Omega_{21}^\dagger \hat k^{-1} \hat \Omega_{21} - 2c\lambda \hat H - \hat S_1\right]\delta y = 0, 
\end{align}
where $\lambda= \gamma + i\omega$ is the eigenvalue that determines stability of the Stokes wave. The quasiperiodic case amounts to putting the subscript $\mu$ onto the respective operators:
\begin{align}
    \left[ \lambda^2 \hat \Omega_{21,\mu}^\dagger \hat k^{-1}_{\mu} \hat \Omega_{21,\mu} - 2c\lambda \hat H_\mu - \hat S_{1,\mu}\right]p(u) = 0,\label{quasiQEP}
\end{align}
where we note that $\hat \Omega_{21,\mu}^\dagger \hat k^{-1}_{\mu} \hat \Omega_{21,\mu}$, $\hat S_{1,\mu}$ are self-adjoint operators and $\hat H_{\mu}$ is skew-adjoint. One can show that the dispersion relation for linear waves is,
\begin{align}
    \omega_\mu(k) = c (k+\mu) \pm \sqrt{g |k + \mu|}.
\end{align}

For numerical purposes it is more convenient to work with the first order system, given by:
\begin{align}
\hat k_\mu \delta \phi  &= \lambda \hat \Omega_{21,\mu} \delta y \\
\hat S_{1,\mu} \delta y &= \lambda \hat \Omega_{21,\mu}^\dagger \delta \phi - 2c\lambda \hat H_\mu \delta y,
\end{align}
or written in the block operator form:
\begin{align}
    \left(\hat M_{\mu} - \lambda \hat L_{\mu}\right)\left[ 
    \begin{array}{c}
         \delta \phi  \\
         \delta y 
    \end{array}
    \right] = 0,
\end{align}
and we proceed following the result in~\cite{dyachenko_semenova2022} to form the eigenvalue problem, however we include the FFH method~\cite{deconinck2006computing}.
We write the inverse operator to $\hat L_{\mu}$ as 
follows:
\begin{align}
    \hat L_\mu =
\left[
\begin{array}{c|c}
0  & \hat \Omega_{21,\mu} \\ \hline
\hat \Omega_{21,\mu}^\dagger  & -2c\hat H_\mu
\end{array}
\right]
\,\,\mbox{and}\,\,
\hat L_{\mu}^{-1}  =
\left[
\begin{array}{c|c}
2c\hat R^\dagger_{12,\mu}\hat H_\mu \hat R_{12,\mu}  & \hat R_{12,\mu}^\dagger \\ \hline
\hat R_{12,\mu}  & 0
\end{array}
\right],
\label{OperLinMu}
\end{align}
and note that $\hat L^{-1}_\mu = \hat R^\dagger_{\mu} J^{-1}_{\mu} \hat R_{\mu}$, where we defined the operators from~\eqref{OperLin} with Floquet multiplier $\mu$. We introduce the change of basis 
matrix $\hat R_{\mu}$ and the vector of canonical coordinates, ${\bf w} = \left(\delta\mathcal{P},\delta y \right)^T = \hat \Omega_{21}^\dagger {\bf v}$, satisfying ${\bf v} = \hat R_{\mu}^\dagger {\bf w}$ where ${\bf v} = (\delta \phi, \delta y)^T$, by the relations:
\begin{align}
    \hat R_{\mu} = \hat \Omega_{\mu}^{-1} = \diag \left( \hat R_{12,\mu} 1\right), \quad \mbox{and} \quad \hat R^{-1}_{\mu} = \hat \Omega_{\mu} = \diag \left( \hat \Omega_{21,\mu}, 1\right),
\end{align}
and 
\begin{align}
    \hat J_{\mu} =  \left[
\begin{array}{c|c}
0  & 1 \\ \hline
1  & -2c\hat H_{\mu}
\end{array}
\right].
\end{align}
We write the shift--invert eigenvalue problem for the linearization on canonical coordinates ${\bf w}$ with purely imaginary shift $i\sigma$ as follows:
\begin{align}
\left(\hat R_{\mu} \hat M_{\mu} \hat R_{\mu}^\dagger - i\sigma J_\mu\right)^{-1} \hat J_\mu \,{\bf w} = \frac{1}{\lambda - i\sigma}\, {\bf w}.
\end{align}
The inverse operator in the left hand side is applied to an
arbitrary $2\pi$-periodic function via the formula:
\begin{align}
    \left(\hat R_{\mu} \hat M_{\mu} \hat R_{\mu}^\dagger - i\sigma J_\mu\right)^{-1} 
    \left[
    \begin{array}{c}
         f  \\
         g 
    \end{array}
    \right] = 
    \left[
    \begin{array}{c}
         \hat Q_{\mu} f \\
         0 
    \end{array}
    \right] + 
    \left[
    \begin{array}{c}
         i\sigma \hat Q_{\mu}  \\
         1 
    \end{array}
    \right]\hat S_{2,\mu}^{-1} \left(g + i\sigma \hat Q_{\mu}f\right).
\end{align}
Here $\hat Q_\mu := \hat \Omega_{21,\mu}^\dagger \hat k_{\mu}^{-1} \hat \Omega_{21,\mu}$ and 
$\hat S_{2,\mu} := \hat S_{1,\mu} + 2ic\sigma\hat H_{\mu} + \sigma^2\hat Q_{\mu}$ are applied in 
$O(N\log N)$ flops. The inverse to $\hat S_{2,\mu}$ is 
found from the solution of the linear operator equation by using 
the minimal residual (MINRES) method~\cite{saad1992numerical} with a  diagonal preconditioner in Fourier space. Applying the inverse requires $O(mN\log N)$ flops, where $m\ll N$ is the number of iterations in MINRES.

\section{Conclusions}
We considered two eigenvalue problems 
that are fundamental to the study of Stokes waves and their stability.  Numerical schemes are developed that allow to solve them efficiently via a matrix-free pseudospectral method with $O(N\log N)$ numerical complexity and suitable for the study of quasiperiodic eigenfunctions.

The first eigenvalue problem that we considered is posed for the linearization of the Babenko equation.
The eigenfunctions of the linearized Babenko operator can be classified in two types -- even and odd. Both types of eigenfunctions originate from $\cos(kx)$ and $\sin(kx)$ in a flat surface associated with a double eigenvalue (except for the $k=0$ case). In a finite amplitude Stokes wave, the double eigenvalue breaks into a pair of simple eigenvalues. Each one is associated with an even and an odd eigenfunction (shown in the left and right panels of Figs.~\ref{fig:eigS1}--\ref{fig:eigfS1}). For all waves, zero is an eigenvalue associated with an odd eigenfunction, $y_u$, which is the derivative of the Stokes wave solution. 
In the waves at bifurcation points of $\hat S_1$, the zero eigenvalue has double multiplicity and corresponds to one of the eigenvalues cross zero. 
Only eigenvalues associated with even eigenfunctions cross zero. The bifurcations occur at Stokes waves either at a turning-point in the speed (for co-periodic eigenfunctions), 
or at a secondary bifurcation to quasiperiodic branch (e.g. period doubling).

The second eigenvalue problem that was considered is posed for the linearization of the dynamical equations around a Stokes wave and governs its stability. It is shown how the stability problem can be interpreted as a second order pseudo-differential PDE and solved numerically by means of canonical conformal variables method~\cite{dyachenko_semenova2022} with the extension to quasiperiodic disturbances via the FFH approach~\cite{deconinck2006computing}.
We considered stability spectrum for waves studied in previous Refs.~\cite{korotkevich2022superharmonic,dyachenko_semenova2022}  with the FFH approach and have shown accuracy and performance of the new method. We extended the results of the aforementioned works and reported novel physically relevant findings in the Ref.~\cite{deconinck2022instability}.
Additional results for dynamical spectrum in the vicinity of the limiting wave for the Benjamin-Feir, localized and high-frequency branches are a subject of ongoing research.  

\section{Acknowledgements}
We would like to thank Eleanor Byrnes,  Prof. Thomas Bridges, Prof. Bernard Deconinck,  Prof. Alexander Korotkevich and Prof. Pavel Lushnikov for fruitful discussion and valuable suggestions. S.D would like to thank the Newton Institute at Cambridge University
for hospitality during the ``Dispersive hydrodynamics: mathematics, simulation and experiments'' program.
S.D. thanks the Department of Applied Mathematics at the U. of Washington for hospitality. A.S and S.D. acknowledge the FFTW project and its authors~\cite{frigo2005design} as well as the entire GNU Project. 
A.S. thanks the Institute for Computational and Experimental Research in Mathematics, Providence, RI, being resident during the ``Hamiltonian Methods in Dispersive and Wave Evolution Equations" program supported by NSF-DMS-1929284.



\bibliographystyle{elsarticle-num} 
\bibliography{StokesWave}

\end{document}